\documentclass[11pt,leqno,twoside]{article}

\usepackage{amsfonts,amsmath,amsthm,amssymb}
\usepackage{enumerate}
\usepackage{graphics,graphicx,subfigure}

\usepackage{color}
\usepackage{todonotes}
\usepackage{cancel}
\usepackage{url}
\usepackage{hyperref}
\usepackage{makeidx}
\usepackage{showidx}
\usepackage{multicol}        
\usepackage{xspace}
\usepackage{stmaryrd}        
\usepackage{pifont}          
\usepackage{fancybox}        
\usepackage{bm}

\usepackage{algorithm}
\usepackage{algorithmic}

 \usepackage{fancyhdr}
\theoremstyle{plain}
 \theoremstyle{definition}
 \newtheorem{lem}{Lemma}
 \newtheorem{defn}[lem]{Definition}
 \newtheorem{thm}[lem]{Theorem}
 \newtheorem{prop}[lem]{Proposition}
 \newtheorem{cor}[lem]{Corollary}
 \newtheorem{notn}[lem]{Notations}
 \newtheorem{pb}[lem]{Problem}
 \newtheorem{form}[lem]{Formulae}
 
 \newtheorem{rk}[lem]{Remark}
 \newtheorem*{com}{Comment}
 \newtheorem*{ex}{Example}
 \theoremstyle{remark}
 \newtheorem{quest}{{\bf Question}}

 \newcommand{\blem}{\begin{lem}}
 \newcommand{\elem}{\end{lem}}
 \newcommand{\bdefn}{\begin{defn}}
 \newcommand{\edefn}{\end{defn}}
 \newcommand{\bthm}{\begin{thm} }
 \newcommand{\ethm}{\end{thm}}
 \newcommand{\bprop}{\begin{prop}}
 \newcommand{\eprop}{\end{prop}}
 \newcommand{\bcor}{\begin{cor}}
 \newcommand{\ecor}{\end{cor}}
 \newcommand{\bnotn}{\begin{notn}}
 \newcommand{\enotn}{\end{notn}}
 \newcommand{\bpb}{\begin{pb}}
 \newcommand{\epb}{\end{pb}}
 \newcommand{\bform}{\begin{form}}
 \newcommand{\eform}{\end{form}}
 \newcommand{\brk}{\begin{rk}}
 \newcommand{\erk}{\end{rk}}
 \newcommand{\bcom}{\begin{com}}
 \newcommand{\ecom}{\end{com}}
 \newcommand{\bex}{\begin{ex}}
 \newcommand{\eex}{\end{ex}}
 \newcommand{\bpf}{\begin{proof}}
 \newcommand{\epf}{\end{proof}}
 \newcommand{\bquest}{\begin{quest} \em}
 \newcommand{\equest}{\end{quest}}

\newcommand{\va}{{\bf a}}

\newcommand{\vc}{{\bf c}}

\newcommand{\ve}{{\bf e}}

\newcommand{\vu}{{\bf u}}
\newcommand{\vv}{{\bf v}}

\newcommand{\vx}{{\bf x}}
\newcommand{\vy}{{\bf y}}

\newcommand{\vA}{{\bf A}}
\newcommand{\vB}{{\bf B}}

\newcommand{\vI}{{\bf I}}

\newcommand{\vM}{{\bf M}}

\newcommand{\vW}{{\bf W}}
\newcommand{\vX}{{\bf X}}
\newcommand{\vY}{{\bf Y}}
\newcommand{\vZ}{{\bf Z}}

\newcommand{\cA}{\mathcal{A}}
\newcommand{\cB}{\mathcal{B}}

\newcommand{\cN}{\mathcal{N}}

\newcommand{\cV}{\mathcal{V}}

\newcommand{\bR}{\mathbb{R}}

\newcommand{\be}{\begin{equation}}
\newcommand{\ee}{\end{equation}}
\newcommand{\bal}{\begin{align}}
\newcommand{\eal}{\end{align}}
\newcommand{\ba}{\begin{align*}}
\newcommand{\ea}{\end{align*}}
\newcommand{\bmx}{\begin{matrix}}
\newcommand{\emx}{\end{matrix}}
\newcommand{\bbmx}{\begin{bmatrix}}
\newcommand{\ebmx}{\end{bmatrix}}
\newcommand{\bpmx}{\begin{pmatrix}}
\newcommand{\epmx}{\end{pmatrix}}
\newcommand{\bvmx}{\begin{vmatrix}}
\newcommand{\evmx}{\end{vmatrix}}

\newcommand{\ol}{\overline}

\newcommand{\f}{\frac}

\newcommand{\imp}{\Longrightarrow}

\newcommand{\inc}{\subseteq}

\newcommand{\tto}{\longrightarrow}

\newcommand{\lbt}{\llbracket}
\newcommand{\rbt}{\rrbracket}

\newcommand{\sgn}{\mathrm{sgn}}

\newcommand{\argmin}{{\rm argmin}\,}
\newcommand{\argmax}{{\rm argmax}\,}

\newcommand{\maximize}{{\rm maximize}\,}

\title{\vspace{-25mm}Jointly Low-Rank and Bisparse Recovery:
\\
Questions and Partial Answers\medskip\hrule height 1.2pt \vspace{-6mm}}
\author{Simon Foucart\footnote{S. F. partially supported by NSF grants
    DMS-1622134 and DMS-1664803. L. J. is funded by the FNRS, Belgium.},
    \, R\'emi Gribonval, Laurent Jacques, and Holger Rauhut}
\date{\vspace{-8mm}\rule{100mm}{0.8pt}}

\newcommand\shorttitle{Jointly Low-Rank and Bisparse Recovery}
\newcommand\authors{S. Foucart, R. Gribonval, L. Jacques, H. Rauhut}

\begin{document}
\maketitle

\vspace{-15mm}
\begin{abstract}
We investigate the problem of recovering jointly $r$-rank and $s$-bisparse matrices from as few linear measurements as possible, considering arbitrary measurements as well as rank-one measurements. In both cases, we show that $m \asymp r s \ln(en/s)$ measurements make the recovery possible in theory, meaning via a nonpractical algorithm.

In case of arbitrary measurements, we investigate the possibility of achieving practical recovery via an iterative-hard-thresholding algorithm when $m \asymp r s^\gamma \ln(en/s)$ for some exponent $\gamma > 0$. We show that this is feasible for $\gamma = 2$, and that the proposed analysis cannot cover the case $\gamma \leq 1$. The precise value of the optimal exponent $\gamma \in [1,2]$ is the object of a question, raised but unresolved in this paper, about head projections for the jointly low-rank and bisparse structure.

Some related questions are partially answered in passing.
For rank-one measurements,
we suggest on arcane grounds an iterative-hard-thresholding algorithm modified to exploit the nonstandard restricted isometry property obeyed by this type of measurements.
\end{abstract}\vspace{-6mm}

\tableofcontents

\linespread{1}

\section{Introduction}

This whole article is concerned with the inquiry below.

{\em {\bf Main Question.}
What is the minimal number of linear measurements needed to recover jointly $r$-rank and $s$-bisparse symmetric $n \times n$ matrices via an efficient algorithm}?

This minimal number of measurements will be called sample complexity.
We will show that it is of the order $r s \ln(en/s)$.
Nevertheless, we do not consider the question fully resolved because
of the lack of efficient algorithms for arbitrary measurements and of the limitation of an efficient algorithm to factorized measurements, and thus to the only applications that could support such a structured sensing. 
Settling the question by providing an efficient algorithm applicable to any type of measurements is therefore still open.
Before diving into our investigations, let us start by clarifying a few points.

\vspace{-2mm}
$\bullet$\; What are `jointly $r$-rank and $s$-bisparse symmetric $n \times n$ matrices'?\\
In this article, we consider exclusively matrices $\vX \in \bR^{n \times n}$ that are symmetric,
i.e., $\vX^\top = \vX$.
The set of $r$-rank (symmetric) matrices will be denoted as
\be
\Sigma^{[r]} := \left\{
\vX \in \bR^{n \times n}: \vX^\top = \vX, \; {\rm rank}(\vX) \le r
\right\}
\ee
and the set of $s$-bisparse (symmetric) matrices will be denoted as
\be
\Sigma_{(s)} := \left\{
\vX \in \bR^{n \times n}: \vX^\top = \vX, \; \vX_{\ol{S \times S}} = {\bf 0}
\mbox{ for some $S \inc \lbt 1:n \rbt$ with $|S|=s$}
\right\},
\ee
where $\vM_\Omega = {\bf 0}$ for $\vM \in \bR^{n \times n}$ and $\Omega \inc \lbt 1:n \rbt \times \lbt 1:n \rbt$ means that all entries of $\vM$ indexed by $\Omega$ are~zeros, and $\overline{\Omega}$ stands for the complement of $\Omega$.

Hence, the jointly $r$-rank and $s$-bisparse (symmetric) matrices we are interested in are elements of
\be
\Sigma_{(s)}^{[r]} := \Sigma^{[r]} \cap \Sigma_{(s)}.
\ee
We will often use the fact that $\Sigma_{(s)}^{[r]} + \Sigma_{(s)}^{[r]} \inc \Sigma_{(2s)}^{[2r]}$.

Note that, as described below, $\Sigma^{[1]}_{(s)}$ is for instance the set associated with the lifting of sparse signals to rank-one matrices when one is interested in their recovery from phaseless (complex) measurements \cite{IwenVW}, while for $r>1$, any matrix of $\Sigma^{[r]}_{(s)}$ describes a quadratic function of both few variables and few quadratic terms whose sampling and recovery --- an important problem in, e.g., approximation theory and high-dimensional statistics --- are related to the Main Question \cite{Fou19,DPW11}.

\vspace{-2mm}
$\bullet$\; What are the `linear measurements' considered?\\
They can be of the arbitrary type
\be 
\label{ArbMst}
y_i = \langle \vX, \vA_i \rangle_F = {\rm tr} (\vA_i^\top \vX),
\qquad i \in \lbt 1:m \rbt,
\ee
or of the specific (rank-one) type
\be
\label{SpeMst}
y_i = \langle \vX \va_i, \va_i \rangle = {\rm tr} (\va_i \va_i^\top \vX),
\qquad i \in \lbt 1:m \rbt .
\ee
Generically, we write $\vy = \bm{\cA}(\vX)$,
where $\bm{\cA}: \bR^{n \times n} \to \bR^m$ is a linear map.

\vspace{-2mm}
$\bullet$\; What is meant by `recover'?\\
More than just finding a map $\Delta: \bR^m \to \bR^{n \times n}$ such that $\Delta(\bm{\cA}(\vX)) = \vX$ for all $\vX \in \Sigma_{(s)}^{[r]}$.
Indeed, we require the recovery procedure to be stable and robust, in the sense that we want
\be
\label{ObjStaRob}
\|\vX - \Delta(\bm{\cA}(\vX) + \ve)\|
\le C \min_{\vZ \in  \Sigma_{(s)}^{[r]}} \|\vX - \vZ \| + D \|\ve\| 
\ee
to hold for all $\vX \in \bR^{n \times n}$ and all $\ve \in \bR^m$.
We give ourselves some freedom on the choice of the three norms appearing in \eqref{ObjStaRob}.
We also require the recovery procedure to be implementable by a practical algorithm, that is, an efficient algorithm 
whose run-time is at most polynomial in $n$~and~$m$ (ideally, a polynomial of low degree, of course).

In our study of the Main Question, we faced the following puzzle.

\bquest
\label{Q2}
Given a positive constant $c \le 1$,
for which value of $s'$, depending on $s$,
can one find a practical algorithm that constructs,
for each symmetric matrix $\vM \in \bR^{n \times n}$,
an index set $S'$ of size~$s'$ such that
\be
\|\vM_{S' \times S'}\|_F^2 \ge c \, \max_{|S|=s}
\|\vM_{S \times S}\|_F^2 \; ?
\ee

\equest

In reality, the relevant question for our goal is broader.
It involves the projection $P^{[r]}$ onto~$\Sigma^{[r]}$.

\bquest
\label{Q3}
Given a positive constant $c \le 1$,
for which value of $s'$, depending on $s$,
can one find a practical algorithm that constructs,
for each symmetric matrix $\vM \in \bR^{n \times n}$,
an index set $S'$ of size~$s'$ such that,
with $r'$ proportional to $r$,
\be
\|P^{[r']}(\vM_{S' \times S'})\|_F^2 \ge c \, \max_{|S|=s}
\|P^{[r]}(\vM_{S \times S})\|_F^2 \; ?
\ee
\equest

If $s'$ could be chosen proportional to $s$ in Question \ref{Q3},
then the Main Question could be answered with $m \asymp r  s \ln(en/s)$ measurements satisfying the so-called restricted isometry property (see below).
This is shown in Section~\ref{SecSCPract}.

We come up with partial answers to the above questions: in Proposition~\ref{PropHeadstos2} we show that for $c=1$ the answer to Question~\ref{Q2} is positive with $s' = s^2$, but that it is negative for any $c>0$ when $s' = O(s)$. Combined with the results of Section~\ref{SecSCPract} this establishes that the answer to the Main Question is positive with $m \asymp rs^\gamma \ln(en/s)$ and $\gamma=2$, using a practical variant of iterative hard thresholding, and that the proposed analysis cannot cover the case $\gamma \leq 1$.

In principle,
we are more interested in the measurements of type \eqref{SpeMst}.
Indeed, in the particular case $r=1$,
the measurements taken on a matrix of the type $\vX = \vx \vx^\top \in \Sigma_{(s)}^{[1]}$ with an $s$-sparse $\vx \in \bR^n$ would read 
\be
y_i = | \langle \va_i, \vx \rangle |^2,
\qquad i \in \lbt 1:m \rbt.
\ee
This is exactly the framework of sparse phaseless recovery
(except that everything should be written in the complex setting).
In this case, the sample complexity is known \cite{IwenVW} to be of the order $m \asymp s \ln(en/s)$,
although it is unclear if this can be achieved with independent Gaussian vectors $\va_1,\ldots,\va_m \in \bR^{n}$.

\brk Similar problems as studied here appear in the context of low-rank tensor recovery where one would
like to project onto the intersection of two or more low rank structures defined by different matricizations.
It is NP-hard to compute exact projections and efficiently computable approximate projections are not yet good enough
to show low-rank tensor recovery results for corresponding iterative hard thresholding guarantees \cite{rascst17}. 
They are also considered in the context of sparse PCA 
from inaccurate and incomplete measurements where the problem of recovering a low-rank matrix with sparse (or compressible) right-singular vectors is analyzed~\cite{For18}. In this work, a multi-penalty approach called A-T-LAS$_{1,2}$ provably reaches local convergence from a reliable, computable initialization. 
Other locally convergent methods applied to the recovery of row-sparse (or column-sparse) and low-rank matrices are the sparse power factorization (SPF) and its subspace-concatenated variant (SCSPF), see \cite{Lee17}. 
While the latter work assumes a high peak-to-average power ratio on the singular vectors of the observed matrix, \cite{Gep19} recently enlarged the class of recoverable matrices by relaxing this constraint.
\erk

\vspace{-5mm}
\section{Theoretical Sample Complexity}\vspace{-3mm}
\label{SecTSC}

Restricted isometry properties have been central in all sorts of structured recovery problems. It is no surprise that another instance of a restricted isometry property plays a key role here, too.
The proof sketch is deferred to the appendix.

\bthm
\label{ThmStandardRIP}
Suppose $\vA_1,\ldots,\vA_m$ are independent random matrices with independent $\cN(0,1/m)$ entries.
Given $\delta > 0$, there exist two values $C,c > 0$ (only depending on $\delta$), such that, with failure probability at most $2 \exp(- c m)$,
\be
\label{GenRIP}
(1-\delta) \|\vZ\|_F^2
\le \|\bm{\cA}(\vZ)\|_2^2
\le (1+\delta) \|\vZ\|_F^2
\qquad
\mbox{for all } \vZ \in \Sigma_{(s)}^{[r]}
\ee
provided $m \ge C r s \ln(en/s)$.
\ethm
\vspace{-2mm}

For the rest of this section,
we place ourselves in the situation where the measurement map $\bm{\cA}$ satisfies the restricted isometry property \eqref{GenRIP},
which can occur as soon as $m$ is of the order $r s \ln(en/s)$.
We can then propose several robust algorithms that recover $\vX \in \Sigma_{(s)}^{[r]}$ from $\vy = \bm{\cA}(\vX) + \ve$.
The first obvious candidate is
\be 
\label{FirstObv}
\Delta(\vy) = \underset{\vZ \in \Sigma_{(s)}^{[r]}}{\argmin} \|\vy - \bm{\cA}(\vZ)\|_2.
\ee
We immediately see that $\|\vy - \bm{\cA}(\Delta(\vy))\|_2 \le \|\vy - \bm{\cA}(\vX)\|_2 = \|\ve\|_2$,
from where it follows that
\be
\label{TSC1}
\| \bm{\cA}(\vX) - \bm{\cA}( \Delta(\vy)  ) \|_2
\le \| \vy - \bm{\cA}( \Delta(\vy)  ) \|_2 + \|\ve\|_2 
\le 2 \|\ve\|_2, 
\ee
and we finally derive that
\be
\label{TSC2}
\|\vX - \Delta(\bm{\cA}(\vX) + \ve) \|_F
\le \f{1}{\sqrt{1-\delta}} \| \bm{\cA}(\vX) - \bm{\cA}( \Delta(\bm{\cA}(\vX) + \ve)  ) \|_2
\le \f{2}{\sqrt{1-\delta}}\|\ve\|_2.
\ee
However, this scheme is not really an appropriate candidate, since producing $\Delta(\vy)$ is NP-hard in general (see below).
 
After a decade or so of $\ell_1$-norm and nuclear norm minimizations,
the next obvious candidate stands out as
\be
\Delta (\vy) = \underset{\vZ \in \bR^{n \times n}}{\argmin} F(\vZ)
\qquad \mbox{subject to} \quad
\|\vy - \bm{\cA}(\vZ) \|_2 \le \|\ve\|_2,
\ee
where $F$ is a convex function promoting the joint low-rank and bisparsity structure.
The negative results from \cite{OJFEH} indicate that reducing the sample complexity below $\min\{ rn, s^2 \ln(en/s) \}$ is unattainable when $F$ is a positive combination of the $\ell_{1}$-norm and nuclear norm.

What about a variant of iterative hard thresholding?
Consider the sequence $(\vX_k)_{k \ge 0}$ defined by
\be
\label{SimpleIHT}
\vX_{k+1} = P_{(s)}^{[r]}(\vX_k + \bm{\cA}^*(\vy - \bm{\cA} (\vX_k))),
\ee
where the adjoint of $\bm{\cA}$ is given by
$$
\bm{\cA}^* : \vu \in \bR^m \mapsto \sum_{i=1}^m u_i \vA_i \in \bR^{n \times n}
$$
and where $P_{(s)}^{[r]} : \bR^{n \times n} \to \Sigma_{(s)}^{[r]}$ denotes the projection onto $\Sigma_{(s)}^{[r]}$, that is, the operator of best approximation from $\Sigma_{(s)}^{[r]}$.
One can show (see Appendix or \cite{Blu}) that if $\Delta(\vy)$ is defined as a cluster point of $(\vX_k)_{k \ge 0}$,  
then 
\be
\| \vX - \Delta(\bm{\cA}(\vX) + \ve) \|_F \le C \|\ve\|_2 
\ee
holds for all $\vX \in \Sigma_{(s)}^{[r]}$ and all $\ve \in \bR^m$.
Here also the issue is that computing $P_{(s)}^{[r]}$ is NP-hard 
(see Section~\ref{SecAppProj}),
which incidentally justifies the NP-hardness of \eqref{FirstObv}
(think of  $\bm{\cA} = \vI$).
What about replacing  $P_{(s)}^{[r]}$ by an operator of near-best approximation from $ \Sigma_{(s)}^{[r]}$, as in, e.g., \cite{GolDav}?
After all, if there is any chance for  \eqref{ObjStaRob} to hold,
then such an operator must exist (think again of $\bm{\cA} = \vI$).
We will in fact construct such an operator in Subsection \ref{SSecJoint}.
But substituting $P_{(s)}^{[r]}$ by such an operator in the proof of Theorem \ref{ThmSimpleIHT} (see Appendix) is not enough to do the trick.

\section{Optimal Sample Complexity with Factorized  Measurements} 

In this section, we show that the optimal sample complexity can be achieved with a practical algorithm in a rather special measurement framework. This framework being restricted to the specific structure of this sensing procedure, the Main Question remains of interest.

We suppose here that matrices $\vX \in \Sigma_{(s)}^{[r]}$ are acquired via measurements in factorized form, namely
\be
\label{eq:fact-sensing}
y_i = \langle \vX, \vB^\top \vA_i \vB \rangle,
\qquad i \in \lbt 1:m \rbt,
\ee
where $\vA_1,\ldots,\vA_m \in \bR^{p \times p}$ allow for low-rank recovery
and $\vB \in \bR^{p \times n}$ allows for sparse recovery.
The recovery algorithm proceeds in two steps, 
which are both practical, i.e., efficiently implementable.\vspace{-5mm}
\begin{itemize}
\item[1.] Compute $\vY^\sharp \in \bR^{p \times p}$ from $\vy \in \bR^m$ as a solution of the nuclear norm minimization
\[
\underset{\vY \in \bR^{p \times p}}{\rm minimize} \; \|\vY\|_* \qquad \mbox{ subject to } \; \langle \vY, \vA_i \rangle_F = y_i, \quad i \in \lbt 1:m \rbt,
\]
or as the output of another low-rank recovery algorithm such as iterative hard thresholding.
\item[2.] Compute $\vX^\sharp \in \bR^{n \times n}$ from $\vY^\sharp$ as the output of the HiHTP algorithm with measurement map $\bm{\cB}: \vZ \in \bR^{n \times n} \mapsto \vB \vZ \vB^\top \in \bR^{p \times p}$.
\end{itemize} 

Although we refer to \cite{roklwuei16,roflkueiwu18} for the exact
formulation of the hierarchically  structured sparsity hard
  thresholding pursuit (HiHTP) algorithm,
a few words about the concept of hierarchical sparsity are in order before we state our result about the two-step recovery procedure above. 
A matrix is said to be $(s,t)$-hierachical sparse (or simply $(s,t)$-sparse)
if at most $s$ of its columns are nonzero and each of these columns possesses at most $t$ nonzero entries. 
Thus, $s$-bisparse matrices are in particular $(s,s)$-sparse.
The HiHTP algorithm essentially relies on the possibility to compute the projection (operator of best approximation) onto $(s,t)$-sparse matrices.
In contrast to the projection onto $s$-bisparse matrices, this is indeed an easy task:
first, select the $t$ largest absolute entries in each column and calculate the resulting $\ell_2$-norm,
then select the $s$ columns with the largest of these $\ell_2$-norms.

\bthm 
\label{ThmFactorized}
Let $\vA_1,\ldots,\vA_m \in \bR^{p \times p}$ be independent standard Gaussian matrices and let $\vB \in \bR^{p \times n}$ be a standard Gaussian matrix independent of $\vA_1,\ldots,\vA_m$. If
\be
p \asymp s \ln(en/s) \qquad \mbox{ and } \qquad m \asymp r p,
\ee
so that $m \asymp rs \ln(en/s)$, 
then the probability that every $\vX \in \Sigma_{(s)}^{[r]}$ is exactly recovered from  
$y_i = \langle \vX, \vB^\top \vA_i \vB \rangle$,
$i \in \lbt 1:m \rbt$,
via the above two-step procedure
is at least $1-2\exp(-cp)$.
\ethm

\bpf
First, notice that the matrix $\vB \vX \vB^\top \in \bR^{p \times p}$ has rank at most~$r$,
since $\vX$ has rank at most~$r$, and that it satisfies
\be 
\langle \vB \vX \vB^\top , \vA_i \rangle_F
 = {\rm tr} (\vA_i^\top \vB \vX \vB^\top) = {\rm tr}(\vB^\top \vA_i^\top \vB \vX ) 
 = \langle \vX, \vB^\top \vA_i \vB \rangle_F
 = y_i,
 \qquad i \in \lbt 1:m \rbt.
\ee
Since $\vA_1,\ldots,\vA_m \in \bR^{p \times p}$ are independent standard Gaussian matrices and $ m \asymp r p$,
it is by now well-known (see, e.g., \cite{CanPla,kakurate16}) that,
with failure probability at most $\exp(-cm)$,
the matrix $\vB \vX \vB^\top$ is recovered via nuclear norm minimization (or another suitable algorithm), so that  $\vY^\sharp = \vB \vX \vB^\top$.

Second, since the matrix $\vX \in \bR^{n \times n}$ is $(s,s)$-sparse and satisfies 
$\bm{\cB}(\vX) = \vB \vX \vB^\top = \vY^\sharp$,
Theorem~1 of \cite{roklwuei16} implies that
the matrix $\vX$ will be exactly recovered via HiHTP as long as the so-called HiRIP of order $(3s,2s)$ holds.
According to Theorem 1 of \cite{roflkueiwu18},
the latter is satisfied when $\vB$ obeys a standard RIP,
and the latter is indeed fulfilled with failure at most $\exp(-cp)$ by the matrix $\vB $ (or rather by a renormalization of it),  
because $\vB \in \bR^{p \times n}$ is a standard Gaussian matrix with $p \asymp s \ln(en/s)$.

All in all, exact recovery of $\vX$ is guaranteed after the two steps with failure probability bounded by $\exp(-cm) + \exp(-cp) \le 2 \exp(-cp)$. 
\epf

\brk
It is possible to extend Theorem \ref{ThmFactorized} beyond the strictly Gaussian setting.
In~particular,
if $\vA_1,\ldots,\vA_m$ take the form $\vA_i = \va_i \va_i^\top$ for some independent standard Gaussian vectors $\va_i \in \bR^p$,
then the first-step recovery of $\vB \vX \vB^\top$ can still be achieved via nuclear norm minimization (see \cite{CaiZha,kakurate16,kurate17}) or by some modified iterative hard thresholding algorithm (see \cite{FouSub}).
Note that the measurements made on $\vX \in \bR^{n \times n}$ are in this case rank-one measurements given by $y_i = \langle \vX \va'_i, \va'_i \rangle$, where $\va'_i := \vB^\top \va_i$.
\erk
 
\brk
Let us mention that sensing strategies similar to \eqref{eq:fact-sensing} have been proposed before for other objects with related structures or for connected problems. 
For instance, 
when estimating $k$-row-sparse and $r$-rank matrices $\vX \in \bR^{n \times n}$  from $m$ ``nested'' measurements $y_i = \langle \vW \vX, \vA_i \rangle$,
\cite{Soh16} showed that RIP conditions  imposed on $\vW \in \bR^{p \times n}$ and on the linear operator associated with $\vA_1,\ldots,\vA_m$
yield a computationally efficient two-stage method that can (nearly) achieve a minimax lower bound from $m \asymp r \max\{p,n\}$ measurements where $p \asymp k \log (n/k)$, i.e., from $m \asymp \max\{ rk \log( n/k), rn\}$. 
A two-stage sensing strategy has been also proposed in \cite{IwenVW} for the sparse phase retrieval problem. 
In this case, the sensing model is factored into a linear operator with robust null space property and a stable phase retrieval matrix
--- the latter allows to recover a compressed form of the sparse vector, 
using e.g. PhaseLift \cite{Can13},
and then the former allows to recover this vector via any compressive sensing algorithm. 
\erk

\vspace{-5mm}
\section{Towards Practical Sample Complexity}\vspace{-3mm}
\label{SecSCPract}

In most scenarios, the measurement map is not of the factorized type considered in  the previous section,
so the two-step procedure cannot even be executed. 
It is therefore still relevant to search for practical recovery algorithms that can be applied
with arbitrary measurement schemes and study the sample complexity using, e.g., Gaussian measurements.
As mentioned at the end of Section~\ref{SecTSC},
a difficulty occurs when one tries to use a near-best approximation operator instead of the best approximation operator $P_{(s)}^{[r]}$ in the iterative hard thresholding algorithm \eqref{SimpleIHT}.
Such a difficulty was also encountered in model-based compressive sensing.
A workaround was found in \cite{HegIndSch}.
As we will see below, our attempt to imitate it prompted Question \ref{Q3}. 

Let us start with the observation that any of the structures $\Sigma_{(s)}$, $\Sigma^{[r]}$, or $\Sigma_{(s)}^{[r]}$ is a union of subspaces, 
which we generically write as
$$
\Sigma = \bigcup_{V \in \cV_{\Sigma}} V.
$$
Then the projection onto $\Sigma$,
i.e., the operator of best approximation from $\Sigma$ with respect to the Frobenius norm, acts on any $\vM \in \bR^{n \times n}$ via
\be
P_{\Sigma}(\vM)  = P_{V(\vM)}(\vM)
\ee
where
\begin{align}
 V(\vM) 
\label{TailProp}& = \underset{V \in \cV_{\Sigma}}{\argmin} \|\vM - P_V(\vM) \|_F^2\\
\label{HeadProp}  & = \underset{V \in \cV_{\Sigma}}{\rm argmax} \; \|P_V(\vM)\|_F^2 ,
\end{align}
and $P_V$ evidently denotes the orthogonal projection onto the subspace $V$.
By analogy with the vector case,
we can think of \eqref{TailProp} as a `tail' property for the projection $P_\Sigma$
and of \eqref{HeadProp} as a `head' property.
We keep this terminology introduced in \cite{HegIndSch} when relaxing the notion of projection.
Precisely, we shall call an operator $T: \bR^{n \times n} \to \Sigma$ a tail projection for $\Sigma$ with constant $C_T \ge 1$ (or near best approximation from $\Sigma$ with constant $C_T$) if
\be
\label{TailCond}
\|\vM - T(\vM) \|_F \le C_T \|\vM - P_{\Sigma}(\vM) \|_F
\qquad \mbox{for all }\vM \in \bR^{n \times n}.
\ee
We may have to relax this notion further by allowing the operator $T$ to map into a bigger set $\Sigma' \supseteq \Sigma$.
Thus, by tail projection for $\Sigma$ into $\Sigma'$ with constant $C_T$, 
we mean an operator $T: \bR^{n \times n} \to \Sigma'$ which satisfies the tail condition \eqref{TailCond}.
Similarly, an operator $H: \bR^{n \times n} \to \Sigma$ is called a head projection for $\Sigma$ with constant $c_H \le 1$ if
\be
\label{HeadCond}
\|H(\vM) \|_F \ge c_H \|P_{\Sigma}(\vM)\|_F
\qquad \mbox{for all }\vM \in \bR^{n \times n}.
\ee
A head projection for $\Sigma$ into $\Sigma' \supseteq \Sigma$ with constant $c_H$ is
an operator $H: \bR^{n \times n} \to \Sigma'$ which satisfies the head condition \eqref{HeadCond}.

At this point, it is worth mentioning (see Appendix) that the (genuine) projection onto $\Sigma_{(s)}^{[r]}$ acts on any $\vM \in \bR^{n \times n}$ via
\be
\label{FormProj}
P_{(s)}^{[r]}(\vM) = P^{[r]}(\vM_{S_\star \times S_\star}),
\qquad \mbox{where }  
S_\star = \underset{|S|=s}{\rm argmax } \; \|P^{[r]}(\vM_{S \times S})\|_F.
\ee
In Section~\ref{SecAppProj}, we will see that 
we can produce a tail projection for $\Sigma_{(s)}^{[r]}$.

The size of $s'$ for which one can produce a head projection for $\Sigma_{(s)}^{[r]}$ into $\Sigma_{(s')}^{[r']}$ with $r'$ proportional to $r$ is exactly the focus of Question~\ref{Q3}. We state and prove below (in the idealized setting where there is no measurement error) that a variant of iterative hard thresholding --- using such a head projection --- allows to perform joint low-rank and bisparse recovery via  from $m \asymp r s' \ln(en/s)$ measurements. This will be interesting if it can be established that $s' \asymp s^\gamma$ with $\gamma <2$ is feasible. Then, for small $r$ (and in particular in the case of sparse phaseless recovery where $r=1$), $m \asymp r s^\gamma \ln(en/s)$ will be of a smaller order than 
both $rn$ --- the sample complexity of rank-$r$ matrices --- and $s^2 \ln(en/s)$. This last bound is associated with enforcing only the matrix bisparse structure, as ensured by combining Theorem \ref{ThmIHTHeadTail} in the case $r=n$ with Proposition \ref{PropHeadstos2} and Theorem~\ref{ThmStandardRIP} (see below). Quite obviously this last context determines that $\gamma=2$ is feasible (as stated in the abstract) since $s^2 \leq r s^2$.

\bthm
\label{ThmIHTHeadTail}
Let $T$ be a tail projection for $\Sigma_{(s)}^{[r]}$ with constant $C_T \ge 1$
and let $H$ be a head projection for $\Sigma_{(2s)}^{[2r]}$ into $\Sigma_{(s')}^{[r']}$ with constant $c_H \le 1$ which additionally takes the form
\be
\label{FormH}
H(\vM) = P^{[r']}(\vM_{S' \times S'})
\qquad \mbox{for some index set $S'$ (depending on $\vM$) of size $s'$}.
\ee
If $(1+C_T)^2  (1 - c_H^2) < 1$ and if the restricted isometry property \eqref{GenRIP} holds on $\Sigma_{(2s+s')}^{[2r+r']}$ with constant $\delta>0$ small enough to have
\be
\rho := (1+C_T)^2  (1 - c_H^2(1-\delta)^2 + 2 \delta(1+\delta)) < 1,
\ee
then any $\vX \in \Sigma_{(s)}^{[r]}$ acquired from $\vy = \bm{\cA}(\vX)$
is recovered as the limit of the sequence $(\vX_k)_{k \ge 0}$ defined by
\be
\vX_{k+1} = T[\vX_k + H(\bm{\cA}^*(\vy - \bm{\cA}(\vX_k)))].
\ee
\ethm

\bpf
We shall prove that, for any $k \ge 0$, 
\be
\label{ObjAnn}
\|\vX - \vX_{k+1}\|_F^2 \le \rho \, \|\vX - \vX_k\|_F^2.
\ee
The tail property guarantees that
\be
\| [\vX_k+H(\bm{\cA}^*(\vy - \bm{\cA}(\vX_k)))] - \vX_{k+1} \|_F \le C_T
\| [\vX_k+H(\bm{\cA}^*(\vy - \bm{\cA}(\vX_k)))] - \vX \|_F 
\ee
and the triangle inequality then yields\footnote{It is probably possible to replace $1+C_T$ by a constant arbitrarily close to $1$ if $T$ mapped into $\Sigma_{(s'')}^{[r'']}$ with $r''$ and~$s''$ proportional to $r$ and $s$ (with proportionality constant increasing when $C_T$ decreases),
as in \cite{ShenLi} for the sparse vector case and in \cite{FouSub} for the low-rank matrix case.
This would allow us to eliminate the condition $(1+C_T)^2  (1 - c_H^2) < 1$.}
\be
\label{StepTI}
\|\vX - \vX_{k+1}\|_F  \le (1+C_T) \| [\vX_k+H(\bm{\cA}^*(\vy - \bm{\cA}(\vX_k)))] - \vX \|_F .
\ee
We now concentrate on bounding $\| [\vX_k+H(\bm{\cA}^*(\vy - \bm{\cA}(\vX_k)))] - \vX \|_F = \| \vZ - H(\bm{\cA}^*\bm{\cA}(\vZ)) \|_F$,
where we have set $\vZ := \vX -\vX_k \in \Sigma_{(2s)}^{[2r]}$.
By expanding the square, we obtain
\begin{align}
\label{ExpSq}
 \| \vZ   - H(\bm{\cA}^*\bm{\cA}(\vZ)) \|_F^2
 & 
 =  \| \vZ\|_F^2 
 + \| H(\bm{\cA}^*\bm{\cA}(\vZ)) \|_F^2
 - 2 \langle \vZ, H( \bm{\cA}^*\bm{\cA}(\vZ) )\rangle_F 
 \\
 \nonumber
 & =  \| \vZ\|_F^2  +  \| H(\bm{\cA}^*\bm{\cA}(\vZ)) \|_F^2
 - 2 \langle \bm{\cA}^* \bm{\cA}(\vZ), H( \bm{\cA}^*\bm{\cA}(\vZ) )\rangle_F\\
 \nonumber
& \phantom{ =  \| \vZ \|_F^2 +  \| H(\bm{\cA}^*\bm{\cA}(\vZ)) \|_F^2}\; - 2 \langle \vZ - \bm{\cA}^* \bm{\cA}(\vZ), H( \bm{\cA}^*\bm{\cA}(\vZ) )\rangle_F.
\end{align}
In view of the form \eqref{FormH} of the head projection, 
followed by the facts that $P^{[r']}$ acts locally as an orthogonal projection
and that it preserves the bisupport of a matrix, 
we observe that
\begin{align}
\label{ObsH}
 \| H(\bm{\cA}^*\bm{\cA}(\vZ)) & \|_F^2
  = \langle  P^{[r']}( \bm{\cA}^*\bm{\cA}(\vZ) _{S' \times S'}),  P^{[r']}( \bm{\cA}^*\bm{\cA}(\vZ) _{S' \times S'}) \rangle_F\\
 \nonumber
 & = \langle   \bm{\cA}^*\bm{\cA}(\vZ) _{S' \times S'},  P^{[r']}( \bm{\cA}^*\bm{\cA}(\vZ) _{S' \times S'}) \rangle_F
  = \langle \bm{\cA}^*\bm{\cA}(\vZ)  ,  P^{[r']}( \bm{\cA}^*\bm{\cA}(\vZ) _{S' \times S'}) \rangle_F\\
 \nonumber
 & = \langle \bm{\cA}^*\bm{\cA}(\vZ)  ,  H(\bm{\cA}^*\bm{\cA}(\vZ))) \rangle_F.
\end{align}
Substituting the latter into \eqref{ExpSq} gives
\be
\label{ExpSq2}
 \| \vZ  - H(\bm{\cA}^*\bm{\cA}(\vZ)) \|_F^2
   =  \| \vZ\|_F^2 
 - \| H(\bm{\cA}^*\bm{\cA}(\vZ)) \|_F^2
  - 2 \langle \vZ - \bm{\cA}^* \bm{\cA}(\vZ), H( \bm{\cA}^*\bm{\cA}(\vZ) )\rangle_F.
\ee
The inner product term is small in absolute value.
Indeed, in view of Lemma \ref{LemSimpleIHT} (see Appendix), we have
\be
\label{IPsmall}
|\langle \vZ - \bm{\cA}^* \bm{\cA}(\vZ), H( \bm{\cA}^*\bm{\cA}(\vZ) )\rangle_F|
 \le \delta \|\vZ\|_F \| H( \bm{\cA}^*\bm{\cA}(\vZ) ) \|_F
 \le \delta (1+\delta) \|\vZ\|_F^2,
\ee
where the bound on $\| H( \bm{\cA}^*\bm{\cA}(\vZ) ) \|_F$ followed from the observation \eqref{ObsH} and the restricted isometry property \eqref{GenRIP},
according to
\begin{align}
\| H( \bm{\cA}^*\bm{\cA}(\vZ) ) \|_F^2
& = \langle \bm{\cA}^*\bm{\cA}(\vZ)  ,  H(\bm{\cA}^*\bm{\cA}(\vZ)) \rangle_F
 = \langle\bm{\cA}(\vZ)  ,  \bm{\cA}( H(\bm{\cA}^*\bm{\cA}(\vZ))) \rangle_F\\
\nonumber
 & \le \|\bm{\cA}(\vZ)\|_F \| \bm{\cA}( H(\bm{\cA}^*\bm{\cA}(\vZ)))  \|_F
 \le (1+\delta) \|\vZ\|_F \| H( \bm{\cA}^*\bm{\cA}(\vZ) ) \|_F.
\end{align}
It now remains to prove that $\| H( \bm{\cA}^*\bm{\cA}(\vZ) ) \|_F^2$ is large,
and this is where the head condition comes into play.
Precisely, assuming that $\vZ$ is supported on $S'' \times S''$ with $|S''| \le 2s$,
we know on the one hand that 
\be
\label{StepHead} 
\| H( \bm{\cA}^*\bm{\cA}(\vZ) ) \|_F \ge c_H 
\| P^{[2r]} ( \bm{\cA}^* \bm{\cA}(\vZ)_{S'' \times S''} ) \|_F.
\ee
On the other hand, using in particular the restricted isometry property \eqref{GenRIP} and Von Neumann's trace inequality combined with the fact that $\vZ$ has rank at most $2r$, we obtain
\begin{align}
\label{StepVN}
(1-\delta) \|\vZ\|_F^2 & \le \| \bm{\cA}(\vZ) \|_2^2
= \langle \vZ, \bm{\cA}^* \bm{\cA} (\vZ) \rangle_F
= \langle \vZ, \bm{\cA}^* \bm{\cA} (\vZ)_{S'' \times S''} \rangle_F\\
\nonumber 
& \le \sum_{i=1}^{2r} \sigma_i( \vZ ) \sigma_i ( \bm{\cA}^* \bm{\cA} (\vZ)_{S'' \times S''} )
\le \left[ \sum_{i=1}^{2r} \sigma_i( \vZ )^2 \right]^{1/2}
\left[ \sum_{i=1}^{2r} \sigma_i ( \bm{\cA}^* \bm{\cA} (\vZ)_{S'' \times S''} )^2 \right]^{1/2}\\
\nonumber
& =  \|\vZ\|_F \|P^{[2r]}(\bm{\cA}^* \bm{\cA} (\vZ)_{S'' \times S''})\|_F.
\end{align}
Combining \eqref{StepHead} and \eqref{StepVN} yields
\be
\label{Hbig}
\| H( \bm{\cA}^*\bm{\cA}(\vZ) ) \|_F \ge c_H (1-\delta) \|\vZ\|_F.
\ee
Substituting \eqref{Hbig} and \eqref{IPsmall} into \eqref{ExpSq2},
we deduce that
\be
\| \vZ  - H(\bm{\cA}^*\bm{\cA}(\vZ)) \|_F^2
\le (1 - c_H^2(1-\delta)^2 + 2 \delta(1+\delta)) \, \|\vZ\|_F^2.
\ee
Finally, using \eqref{StepTI}, we arrive that
\be
\|\vX - \vX_{k+1}\|_F^2 \le (1+C_T)^2  (1 - c_H^2(1-\delta)^2 + 2 \delta(1+\delta)) \, \|\vX - \vX_k\|_F^2,
\ee
which is the objective announced in \eqref{ObjAnn}.
\epf

\section{Tail and Head Projections}
\label{SecAppProj}

In this section, we gather some information about the construction of computable tail and head projections for each of the structures $\Sigma^{[r]}$, $\Sigma_{(s)}$, and $\Sigma_{(s)}^{[r]}$.
We work under the implicit assumption that the domain of all these projections is the space of symmetric matrices,
i.e., the projections are only applied to matrices $\vM \in \bR^{n \times n}$ satisfying $\vM^\top = \vM$.

\subsection{Low-rank structure}

There is no difficulty whatsoever here ---
even the exact projection $P^{[r]}: \bR^{n \times n} \to \Sigma^{[r]}$ is accessible.
Indeed, it is well known that if $\vX \in \bR^{n \times n}$ has singular value decomposition 
\be
\vX = \sum_{i=1}^n \sigma_i(\vX) \vu_i \vv_i^\top
\ee
where the singular values $\sigma_1(\vX) \ge \cdots \ge \sigma_n(\vX) \ge 0$
are arranged in nondecreasing order,
then the projection of $\vX$ onto the set of rank-$r$ matrices is obtained by truncating this decomposition to include only the first $r$ summands,
i.e.,
\be
P^{[r]}(\vX) = \sum_{i=1}^r \sigma_i(\vX) \vu_i \vv_i^\top.
\ee
Note that $P^{[r]}(\vM)$ is symmetric whenever $\vM$ itself is symmetric.

\subsection{Bisparsity structure}
\label{SSecBisparse}

Quickly stated, exact projections for $\Sigma_{(s)}$ are NP-hard,
but there are computable tail projections for $\Sigma_{(s)}$.
Head projections for $\Sigma_{(s)}$ are still NP-hard if they are forced to map exactly into $\Sigma_{(s)}$.

If they are allowed to map into a larger set $\Sigma_{(s')}$,
the situation depends on the order of $s'$ compared to $s$
---  Question~\ref{Q2} in fact asks which value of $s'>s$ allows for a computable head projection.

We provide a few incomplete results related to this situation.

\paragraph{Exact projection.} Finding the exact projection for $\Sigma_{(s)}$ amounts to solving the problem
\be 
\underset{|S|=s}{\maximize} \|\vM_{S \times S}\|_F^2.
\ee
This is NP-hard even with the restriction that $\vM$ is an adjacency matrix of a graph
because it then reduces to the densest $k$-subgraph problem,
which is known to be NP-hard \cite{Manurangsi}.

\paragraph{Tail projections.} 
There is a simple procedure to obtain a practical tail projection for $\Sigma_{(s)}$,
as described below.

\bprop
\label{PropT}
Given a symmetric matrix $\vM \in \bR^{n \times n}$,
let $S_\star$ denote an index set corresponding to $s$ columns of $\vM$ with largest $\ell_2$-norms, i.e.,
\be
S_\star = \underset{|S|=s}{\argmin } \| \vM - \vM_{: \times S} \|_F.
\ee
Then
\be
\| \vM - \vM_{S_\star \times S_\star} \|_F
\le \sqrt{2} \min_{|S|=s} \| \vM - \vM_{S \times S} \|_F.
\ee
\eprop

\bpf
For any index set $T$, the symmetry of $\vM$ imposes that $ \| \vM_{\ol{T} \times T}\|_F^2 = \|\vM_{T \times \ol{T}}\|_F^2$, hence
\be
\|\vM - \vM_{T \times T}\|_F^2 
 = \|\vM_{T \times \ol{T}}\|_F^2  + \|\vM_{\ol{T} \times T}\|_F^2 + \|\vM_{\ol{T} \times \ol{T}}\|_F^2
= 2  \|\vM_{T \times \ol{T}}\|_F^2 + \|\vM_{\ol{T} \times \ol{T}}\|_F^2.
\ee
In view of $\|\vM_{T \times \ol{T}}\|_F^2 + \|\vM_{\ol{T} \times \ol{T}}\|_F^2 = \|\vM_{: \times \ol{T}}\|_F^2 = \| \vM - \vM_{: \times T} \|_F^2$,
we deduce that
\be
\|\vM - \vM_{: \times T}\|_F^2
\le
\|\vM - \vM_{T \times T}\|_F^2
\le 
2 \|\vM - \vM_{: \times T}\|_F^2.
\ee
Applying the latter with $T$ equal to $S_\star$ and with $T$ equal to an arbitrary index set $S$ of size $s$
shows that
\be
\|\vM - \vM_{S_\star \times S_\star}\|_F^2 
\le 2 \|\vM - \vM_{: \times S_\star}\|_F^2 
\le 2 \|\vM - \vM_{: \times S}\|_F^2 
\le 2 \|\vM - \vM_{S \times S}\|_F^2,
\ee
which yields the required result after taking the square root.
\epf

\paragraph{Head projections.}
The literature on the densest $k$-subgraph problem informs us that
finding a head projection for $\Sigma_{(s)}$ is also an NP-hard problem \cite{Manurangsi}.
In our setting, though,
there is room to relax the head projection to map into $\Sigma_{(s')}$ with $s' > s$.
In this regard, Question \ref{Q2} asks if one can actually compute a head projection for $\Sigma_{(s)}$ into $\Sigma_{(s')}$.
We do not have a definite answer for it,
but we prove below that the exponent $\gamma$ in a speculative behavior $s' \asymp s^\gamma$ must lie in $(1,2]$
 --- note that a behavior $s' \asymp s \, {\rm polylog}(s)$ is not excluded.
We then highlight a few observations which feature a nonabsolute constant $c_H$ when $s' \asymp s$.

\begin{algorithm}[t]

\renewcommand{\algorithmicrequire}{\textbf{Input:}}
\renewcommand{\algorithmicensure}{\textbf{Output:}}
\renewcommand{\algorithmicendfor}{\textbf{end}}
\caption{
A head projection $H$ for $\Sigma_{(s)}$ to $\Sigma_{(s^2)}$ with $c_H=1$
\label{alg:head-s-to-2s}
}
\begin{algorithmic}
\REQUIRE A symmetric matrix $\vM \in \bR^{n \times n}$, a sparsity level $s \in \lbt 1:n \rbt$.\vspace{2mm}
\FOR{$i \in \lbt 1:n \rbt$}\vspace{1mm} 
\STATE $C_i := \underset{|C|=s-1,  C \not\ni i}{\argmax} \; \|\vM_{ i \times (\{i\} \cup C)} \|_2$\vspace{1mm}
\STATE $c_i := \|\vM_{ i \times (\{i\} \cup C_i)} \|_2$\vspace{1mm}
\ENDFOR\vspace{2mm}
\STATE $R  := \underset{|R'|=s}{\argmax} \|\vc_{R'}\|$, with $\vc = (c_1,\cdots,c_n)^\top$
\STATE $S' := R \cup (\cup_{i \in R} C_i)$\vspace{2mm}
\RETURN $H(\vM) := \vM_{S' \times S'} \in \Sigma_{(s^2)}$
\end{algorithmic}

\end{algorithm}

\bprop
\label{PropHeadstos2}
The practical algorithm Algorithm~\ref{alg:head-s-to-2s} yields a head projection for $\Sigma_{(s)}$ into $\Sigma_{(s^2)}$ with constant $c_H = 1$.
However, there is no practical algorithm that yields a head projection for $\Sigma_{(s)}$ into $\Sigma_{(s')}$ with absolute constant $c_H >0$ when $s' =O(s)$.
\eprop

\bpf
From the definition of the index sets $\{C_i: 1\leq i \leq n\}$, $R$ and $S'$ in Algorithm~\ref{alg:head-s-to-2s}, for
any index set $S$ with $|S|=s$,
we have
\begin{align}
\|\vM_{S \times S} \|_F^2 
& = \sum_{i \in S} \|\vM_{i \times S} \|_2^2
\le \sum_{i \in S} \|\vM_{ i \times (\{i\} \cup C_i)} \|_2^2
\le \sum_{i \in R} \|\vM_{ i \times (\{i\} \cup C_i)} \|_2^2\\
\nonumber
& = \|\vM_{R \times (R \bigcup \cup_{i \in R} C_i)}  \|_F^2
\le \|\vM_{S' \times S'} \|_F^2,
\end{align}
where the index set $S' = R \bigcup \cup_{i \in R} C_i$ has size at most $s + s(s-1) = s^2$. This proves the first part of the statement.

For the second part of the statement, 
we shall show that if we could compute,
for each $\vM \in \bR^{n \times n}$,
an index set $S'$ with $|S'| \leq Cs$ such that
\be
\label{IfHP}
\|\vM_{S' \times S'} \|_F^2 \ge c_H^2 \max_{|S| = s} \|\vM_{S \times S} \|_F^2,
\ee
then a practical algorithm that yields a head approximation for $\Sigma_{(s)}$ into $\Sigma_{(s)}$ itself would follow,
contradiction the NP-hardness of the latter task.
So let us assume that we have a computable procedure to construct an index set $S'$ as above. 
Looking without loss of generality at the case where $s$ is even and $|S'|=Cs$,
we consider an index set $R \inc S'$ of size $s/2$ corresponding to $s/2$ largest values of $\|\vM_{i \times S'}\|_2$.
By comparing averages, we see that
\be
\label{Av1}
\f{1}{s/2} \|\vM_{R \times S'} \|_F^2 \ge \f{1}{Cs} \|\vM_{S' \times S'} \|_F^2,
\qquad \mbox{i.e.,} \quad
\|\vM_{R \times S'} \|_F^2 \ge \f{1}{2C}\|\vM_{S' \times S'} \|_F^2.
\ee
Next, we consider an index set $C \inc S'$ of size $s/2$ corresponding to $s/2$ largest values of $\|\vM_{R \times j}\|_2$.
 By comparing averages again, we see that
 \be
\label{Av2}
\f{1}{s/2} \|\vM_{R \times C} \|_F^2 \ge \f{1}{Cs} \|\vM_{R \times S'} \|_F^2,
\qquad \mbox{i.e.,} \quad
\|\vM_{R \times C} \|_F^2 \ge \f{1}{2C}\|\vM_{R \times S'} \|_F^2.
\ee
Combining \eqref{Av2}, \eqref{Av1}, and \eqref{IfHP}, we arrive at 
\be
\|\vM_{R \times C} \|_F^2 \ge \f{c_H^2}{4 C^2} \max_{|S| = s} \|\vM_{S \times S} \|_F^2.
\ee
With $T:= R \cup C$, which has size at most $s$, this immediately implies that
\be
\|\vM_{T \times T} \|_F^2 \ge \f{c_H^2}{4 C^2} \max_{|S| = s} \|\vM_{S \times S} \|_F^2,
\ee
meaning that a head approximation for $\Sigma_{(s)}$ into $\Sigma_{(s)}$ can be produced in a practical way.
Since this is not possible,
the second part of the statement is proved.
\epf

Now that we have established the impracticability of head approximations for $\Sigma_{(s)}$ into $\Sigma_{(Cs)}$ with an absolute constant $c_H$,
we examine what can be done when $c_H$ can depend on specific parameters. 

\bprop
\label{PropH1}

Given a symmetric matrix $\vM \in \bR^{n \times n}$, we consider the practical algorithm that returns the matrix $\vM_{T \times T}$ for a set $T := R \cup C$ defined by the union of the index sets of size~$s$
\begin{align}
R & = \underset{|S|=s}{\rm argmax} \; \|\vM_{S \times :}\|_F^2,\\
C & = \underset{|S|=s}{\rm argmax} \; \|\vM_{R \times S}\|_F^2.
\end{align}
This algorithm yields a head projection for $\Sigma_{(s)}$ into $\Sigma_{(2s)}$ with constant $c_H = \sqrt{s/n}$.
\eprop

\bpf
From the definition of $R$ and $C$, it is painless to see that,
for an arbitrary index set $S$ of size~$s$,
\be
\|\vM_{T \times T}\|_F^2 \ge \|\vM_{R \times C}\|_F^2
\ge \f{s}{n} \|\vM_{R \times :}\|_F^2
\ge \f{s}{n} \|\vM_{S \times :}\|_F^2
\ge \f{s}{n} \|\vM_{S \times S}\|_F^2,
\ee
which concludes the proof.
\epf

When $n > s^2$ (which is the most realistic situation from our perspective),
the previous observation is superseded by the following one.

\begin{algorithm}[t]

\renewcommand{\algorithmicrequire}{\textbf{Input:}}
\renewcommand{\algorithmicensure}{\textbf{Output:}}
\renewcommand{\algorithmicendfor}{\textbf{end}}
\caption{
A head projection $H$ for $\Sigma_{(s)}$ with $c_H = 1/\sqrt{s}$
\label{alg:head-s-var-cst}
}
\begin{algorithmic}
\REQUIRE A symmetric matrix $\vM \in \bR^{n \times n}$, a sparsity level $s \in \lbt 1:n \rbt$.\vspace{2mm}
\FOR{$j \in \lbt 1:n \rbt$}\vspace{1mm} 
\STATE $S_j := \underset{|S|=s,  S \ni j}{\argmax} \; \|\vM_{S \times j}\|^2_2$\vspace{1mm}
\ENDFOR\vspace{1mm}
\STATE $j_\star  := \underset{j \in \lbt 1:n \rbt}{\argmax} \|\vM_{S_j \times j}\|_2^2$\vspace{1mm}
\RETURN $H(\vM) := \vM_{S_{j_\star} \times S_{j_\star}} \in \Sigma_{(s)}$
\end{algorithmic}

\end{algorithm}

\bprop
\label{PropH2}
The practical algorithm Algorithm~\ref{alg:head-s-var-cst} yields a head projection for $\Sigma_{(s)}$ with constant $c_H = 1/\sqrt{s}$.
\eprop

\bpf
It is painless to see that, given the definition of Algorithm~\ref{alg:head-s-var-cst}, for an arbitrary index set~$S$ of size $s$,
\be
\|\vM_{S \times S}\|_F^2 
= \sum_{j \in S} \| \vM_{S \times j} \|_2^2
\le \sum_{j \in S} \| \vM_{S_j \times j} \|_2^2
\le  s \| \vM_{S_{j_\star} \times j_\star} \|_2^2
\le s \| \vM_{S_{j_\star} \times S_{j_\star}} \|_F^2,
\ee
which concludes the proof.
\epf

As a final remark,
we show that head projections can be computed for specific symmetric matrices,
e.g., matrices of rank one.

\bprop

Given a symmetric matrix $\vM = \sum_{k=1}^r \vv_k \vv_k^\top \in \bR^{n \times n}$ of rank-$r$, we consider the practical algorithm that returns the matrix $\vM_{S_\star \times S_\star}$,
with $S_\star := S_1 \cup \cdots \cup S_r$, and $S_k$ the index set of $s$ largest absolute entries of $\vv_k$, $1 \leq k \leq r$. This algorithm yields a head projection for $\Sigma_{(s)}$ into $\Sigma_{(rs)}$ with constant $c_H = 1/\sqrt{r}$
when applied to $r$-rank positive semidefinite matrices. 
\eprop

\bpf
Given the definition of $S_\star$, we are going to show that, for any index set $S$ of size $s$,
\be
\label{pros}
\|\vM_{S_\star \times S_\star}\|_F \ge \f{1}{\sqrt{r}} \|\vM_{S \times S}\|_F .
\ee
To do so, we start by writing
\be
M_{i,j}^2 = \left( \sum_{k=1}^r (\vv_k)_i (\vv_k)_j \right)^2
= \sum_{k,\ell=1}^r (\vv_k)_i (\vv_k)_j (\vv_\ell)_i (\vv_\ell)_j. 
\ee
Then, for any index set $T$, in view of
\begin{align}
\|\vM_{T \times T}\|_F^2 & = \sum_{i,j \in T} \sum_{k,\ell=1}^r (\vv_k)_i (\vv_k)_j (\vv_\ell)_i (\vv_\ell)_j
= \sum_{k,\ell=1}^r \sum_{i,j \in T} (\vv_k)_i (\vv_\ell)_i (\vv_k)_j  (\vv_\ell)_j\\
\nonumber
& =\sum_{k,\ell=1}^r  \left( \sum_{i \in T} (\vv_k)_i (\vv_\ell)_i \right)^2,
\end{align}
we derive on the one hand that
\be
\label{OT1H}
\|\vM_{T \times T}\|_F^2
\ge \sum_{k=1}^r  \left( \sum_{i \in T} (\vv_k)_i^2 \right)^2
\ee
and on the other hand,
by the Cauchy--Schwarz inequality applied twice,
that
\be
\label{OT2H}
\|\vM_{T \times T}\|_F^2
\le \sum_{k,\ell=1}^r  \left( \sum_{i \in T} (\vv_k)_i^2 \right) 
\left( \sum_{i \in T} (\vv_\ell)_i^2 \right)
= \left( \sum_{k=1}^r \sum_{i \in T} (\vv_k)_i^2  \right)^2
\le r  \sum_{k=1}^r  \left( \sum_{i \in T} (\vv_k)_i^2 \right)^2.
\ee
Applying \eqref{OT2H} with $T=S$ and using the defining property of each $S_k$ and of $S_\star$, we obtain
\be
\|\vM_{S\times S}\|_F^2 
\le r  \sum_{k=1}^r  \left( \sum_{i \in S} (\vv_k)_i^2 \right)^2 \hspace{-1mm}
\le r  \sum_{k=1}^r  \left( \sum_{i \in S_k} (\vv_k)_i^2 \right)^2 \hspace{-1mm}
\le r  \sum_{k=1}^r  \left( \sum_{i \in S_\star} (\vv_k)_i^2 \right)^2 \hspace{-1mm}
\le r \| \vM_{S_\star \times S_\star} \|_F^2,
\ee
the last inequality being \eqref{OT1H} applied with $T=S_\star$.
The prospective inequality \eqref{pros} is proved.
\epf

\subsection{Joint low-rank and bisparsity structure}
\label{SSecJoint}

Quickly stated, exact projections for $\Sigma_{(s)}^{[r]}$ are NP-hard,
but there are computable tail projections for $\Sigma_{(s)}^{[r]}$.
Head projections for $\Sigma_{(s)}^{[r]}$ are still NP-hard if they are forced to map exactly into $\Sigma_{(s)}^{[r]}$.

If they are allowed to map into a larger set $\Sigma_{(s')}^{[r']}$,
the situation is not settled --- 
this directly relates to Question \ref{Q3}.

We provide a few incomplete results related to this situation.

\paragraph{Exact projections.} 
We already know from Subsection \ref{SSecBisparse} that it is NP-hard to find the exact projection onto $\Sigma_{(s)}^{[r]}$ in general,
since we are talking about exact projection onto $\Sigma_{(s)}$ when $r=n$.
But we are more interested in the case where $r$ is a small constant,
say $r=1$ as a prototype.
Then finding the exact projection onto $\Sigma_{(s)}^{[1]}$ amounts to solving the problem
\be
\underset{|S|=s}{\maximize}\,\|P^{[1]}(\vM_{S \times S})\|_F\ =\ \underset{|S|=s}{\maximize}\, \sigma_{\max}(\vM_{S \times S}).
\ee
Thus, when $\vM$ is a positive semidefinite matrix, we consider the problem  
\be
\underset{\|\vx\|_0 \le s, \|\vx\|_2=1}{\maximize} \langle \vM \vx, \vx \rangle.  
\ee
This is the so-called sparse principal component analysis problem, which is NP-hard \cite{MagdonIsmail}. 

\paragraph{Tail projections.}
There is a fairly simple procedure to create a practical tail projection for~$\Sigma_{(s)}^{[r]}$.
It is based on the availability of tail projections for both $\Sigma^{[r]}$ and $\Sigma_{(s)}$.
The argument is in fact valid for any two `structures' $\Sigma'$ and $\Sigma''$
such that $\Sigma'$ is compatible with a tail projection $T''$ for~$\Sigma''$, in the sense that
\be
\vZ \in \Sigma' \imp T''(\vZ) \in \Sigma'.
\ee
The compatibility applies to the low-rank and bisparsity structures in two different ways:
firstly,
$\Sigma^{[r]}$ is compatible with the tail projection for $\Sigma_{(s)}$ given in Proposition \ref{PropT}, 
by virtue of the fact that a matrix $\vZ$ of rank at most $r$ has all its submatrices $\vZ_{S \times S}$ of rank at most $r$, too;
secondly,
$\Sigma_{(s)}$ is compatible with the exact projection for $\Sigma^{[r]}$,
by virtue of the fact that a matrix $\vZ$ supported on $S \times S$ has all its singular vectors supported on $S$, so that $P^{[r]}(\vZ)$ is supported on $S \times S$, too.
Here is the abstract statement valid for arbitrary structures $\Sigma'$ and $\Sigma''$.

\bprop
Let $T'$ and $T''$ be tail projections for $\Sigma'$ and $\Sigma''$ with constants $C_{T'}$ and $C_{T''}$.
If $\Sigma'$ is compatible with $T''$,
then $T'' \circ T'$ is a tail projection for $\Sigma' \cap \Sigma''$ with constant $C_{T'} + C_{T''} + C_{T'} C_{T''}$.
\eprop

\bpf
We first remark that the compatibility condition ensures that $T'' \circ T'$ maps into $\Sigma' \cap \Sigma''$.
Let $\vM \in \bR^{n \times n}$ and let $P(\vM)$ denote its exact projection for $\Sigma' \cap \Sigma''$.
The tail condition for $T'$ implies that
\be
\label{T'}
\|\vM - T'(\vM) \|_F \le C_{T'} \|\vM - P(\vM)\|_F.
\ee
As a result, we obtain
\be
\label{AaR} 
\|T'(\vM) - P(\vM) \|_F
\le \|T'(\vM) - \vM \|_F + \|\vM - P(\vM) \|_F
\le (C_{T'}+1) \|\vM - P(\vM) \|_F.
\ee
The tail condition for $T''$ combined with \eqref{AaR} yields
\be
\label{T''}
\|T'(\vM) - T''(T'(\vM)) \|_F \le C_{T''}\|T'(\vM) - P(\vM) \|_F
\le  C_{T''} (C_{T'}+1) \|\vM - P(\vM) \|_F.
\ee
Using \eqref{T'} and \eqref{T''}, we derive that
\begin{align}
\|\vM - T''(T'(\vM))\|_F & \le \|\vM - T'(\vM)\|_F+\|T'(\vM) - T''(T'(\vM))\|_F\\
\nonumber
& \le (C_{T'} + C_{T''} (C_{T'}+1)) \|\vM - P(\vM) \|_F,
\end{align}
which proves that $T'' \circ T'$ is a tail projection for $\Sigma' \cap \Sigma''$
with the desired constant. 
\epf

\paragraph{Head projections.}
The literature on the sparse principal component analysis problem informs us that
finding a head projection for $\Sigma_{(s)}^{[r]}$ is still an NP-hard problem \cite[Theorem 2]{MagdonIsmail}.
In our setting, though,
there is room to relax the head projection to map into $\Sigma_{(s')}^{[r']}$ with $r' > r$ and $s' > s$.
In this regard, Question \ref{Q3} asks if one can actually compute a head projection for $\Sigma_{(s)}^{[r]}$ into $\Sigma_{(s')}^{[r']}$ with $r' = Cr$.
We do not have a definite answer for it,
but we highlight an observation featuring a nonabsolute constant $c_H$,
based on what was done for the bisparsity structure.

\bprop
Given a symmetric matrix $\vM \in \bR^{n \times n}$ and $r \le s$, the practical algorithm that yields $P^{[r]}(H(\vM))$ for the operator $H$ defined in Algorithm~\ref{alg:head-s-var-cst} is a head projection for $\Sigma_{(s)}^{[r]}$ with constant $c_H = \sqrt{r}/s$.
\eprop

\bpf
Given a symmetric matrix $\vM \in \bR^{n \times n}$, 
we consider the row (or column) index set $S_\star$ of size $s$ 
supporting the non-zero rows (or columns) of $H(\vM) \in \Sigma_{(s)}$ for the operator $H$ defined in Algorithm~\ref{alg:head-s-var-cst}. By Proposition~\ref{PropH2} for any index set $S$ of size $s$, we have
\be
\|\vM_{S_\star \times S_\star}\|_F^2 \ge \f{1}{s} \|\vM_{S \times S}\|_F^2.
\ee
Then, by noticing that the average of the $r$ largest squared singular values of $\vM_{S_\star \times S_\star}$ is larger than the average of all the squared singular values of $\vM_{S_\star \times S_\star}$,
we derive
\be
\| P^{[r]}(\vM_{S_\star \times S_\star}) \|_F^2
\ge \f{r}{s} \|\vM_{S_\star \times S_\star}\|_F^2
\ge \f{r}{s^2} \|\vM_{S \times S}\|_F^2
\ge \f{r}{s^2} \| P^{[r]}(\vM_{S \times S}) \|_F^2.
\ee
The desired result is now proved.
\epf

A similar argument, based on Proposition \ref{PropH1} instead of Proposition \ref{PropH2}, would yield a head projection for $\Sigma_{(s)}^{[r]}$ into $\Sigma_{(2s)}^{[r]}$ with constant $c_H = \sqrt{r/n}$.

\section{Sample Complexity with Rank-One Measurements}

The specific (rank-one) measurements \eqref{SpeMst} do not result in a measurement map $\bm{\cA} : \bR^{n \times n} \to \bR^m$ obeying the standard restricted isometry property \eqref{GenRIP}.
However, it will satisfy the following version featuring the $\ell_1$-norm as an inner norm.
This was established in \cite{CaiZha} when considering the low-rank structure alone.
The proof sketch is deferred to the appendix.
Note that the rank-one measurements \eqref{SpeMst} also satisfy a version of the null space property ensuring
recovery via nuclear norm minimization, see \cite{kakurate16,kurate17}.

\bthm
\label{ThmSpeRIP}
Suppose $\va_1,\ldots,\va_m \in \bR^m$ are independent vectors 
with independent $\cN(0,1/m)$ entries.
Then, with failure probability at most $2 \exp(-c m)$,
\be
\label{SpeRIP}
\alpha \| \vZ \|_F 
\le \left\|  (\va_i^\top \vZ \va_i)_{i=1}^m \right\|_1
\le \beta \| \vZ \|_F
\qquad
\mbox{for all } \vZ \in \Sigma_{(s)}^{[r]},
\ee
provided $m \ge C r s \ln(en/s)$.
The constants $\beta \ge \alpha >0$ are absolute.	
\ethm

The restricted isometry property \eqref{SpeRIP} already guarantees that the specific-sample complexity --- the theoretical one --- is $m \asymp rs \ln(en/s)$, as expected.
Indeed, given $\vy = \bm{\cA}(\vX) + \ve$ for some $\vX \in \Sigma_{(s)}^{[r]}$,
consider the unpractical recovery scheme
\be
\Delta(\vy) = \underset{\vZ \in \Sigma_{(s)}^{[r]}}{\argmin} \| \vy - \bm{\cA}(\vZ)\|_1.
\ee
In a similar spirit to \eqref{TSC1}-\eqref{TSC2}, we can derive that
\be
\|\vX - \Delta(\bm{\cA}(\vX) + \ve) \|_F \le \f{2}{\alpha}\|\ve\|_1.
\ee

For a practical algorithm scheme,
we have in mind an algorithm belonging to the iterative hard thresholding family.
Namely, we can think of constructing a sequence $(\vX_k)$ of matrices in $\Sigma_{(s')}^{[r']}$
by the recursion\footnote{It is `natural' to include the $\sgn$ operator in order to exploit the restricted isometry property with $\ell_1$ inner norm.}
\be
\label{ModIHT}
\vX_{k+1} =
T \left[ \vX_k + \nu_k H(\bm{\cA}^* \sgn(\vy - \bm{\cA} \vX_k)) \right],
\qquad \nu_k =  \f{\|\vy - \bm{\cA} \vX_k \|_1}{\beta^2}.
\ee
Here, the operators $T : \bR^{n \times n} \to \Sigma_{(s')}^{[r']}$ 
and $H: \bR^{n \times n} \to \Sigma_{(s'')}^{[r'']}$,
depending on parameters $r'$, $s'$, $r''$, and $s''$, may be tail and head projections.
It could also be useful to require the operator $T$ to satisfy the property\footnote{The inequality of \eqref{ProsT} implies that $T$ is a tail projection with $C_T = 1+ \eta(C')$, since 
$$
\|\vM - T(\vM) \|_F \le \|\vM - P_{(s)}^{[r]}(\vM) \|_F + \|P_{(s)}^{[r]}(\vM) - T(\vM) \|_F
\le \|\vM - P_{(s)}^{[r]}(\vM) \|_F + \eta(C') \|P_{(s)}^{[r]}(\vM) - \vM \|_F
= C_T \|\vM - P_{(s)}^{[r]}(\vM) \|_F.
$$
} that,
for all $\vX \in \Sigma_{(s)}^{[r]}$ and all $\vZ \in \bR^{n \times n}$,
\be
\label{ProsT}
\| \vX - T(\vZ) \|_F \le \eta(C) \|\vX - \vZ\|_F
\qquad \mbox{with} \quad
\eta(C')  \underset{C' \to \infty}{\tto} 1.
\ee
With $T = P_{(s')}^{[r']}$,
this inequality seems rather intuitive,
but it needs to be formalized
--- keep in mind, however, that $P_{(s')}^{[r']}$ is not accessible.
When considering the low-rank structure alone,
such an inequality has been established and exploited in \cite{FouSub} to prove that an iterative hard thresholding algorithm of the type \eqref{ModIHT} presents the same recovery guarantees as nuclear norm minimization for recovery from measurements of type \eqref{SpeMst}.
The type of inequality \eqref{ProsT} was first put forward for the sparse vector case in \cite{ShenLi}
and it has been exploited in \cite{FouLec} to propose and analyze an iterative hard thresholding algorithm designed for the case when the standard restricted isometry property fails.

There is an additional property that we could require about the operator $T$.
Namely, 
given a matrix $\vM \in \bR^{n \times n}$,
if $T(\vM)$ is supported on $S \times S$,
then
\be 
\label{ProsT2}
T(\vM) = T(\vM_{S' \times S'})
\qquad \mbox{whenever} \quad
S' \supseteq S.
\ee 
This property is true (see Appendix) for $T = P_{(s')}^{[r']}$, which again is inaccessible.

\section{Appendix: Proofs of Auxiliary Results}

This section collects the detailed arguments for some facts that have been stated but not proved in the narrative. 

\paragraph{Restricted isometry properties.}
First, let us concentrate on Theorem \ref{ThmStandardRIP}
and briefly justify that Gaussian measurements of type \eqref{ArbMst} satisfy the standard restricted isometry property \eqref{GenRIP}.
Without going into details, we simply mention that
the classical proof consisting of a concentration inequality followed by a covering argument works ---
the key being to estimate the covering number of the `ball' of $\Sigma_{(s)}^{[r]}$
essentially as in \cite[Lemma 3.1]{CanPla} with the addition of a union bound.

Next, let us concentrate on Theorem \ref{ThmSpeRIP} and briefly justify that Gaussian rank-one measurements of type \eqref{SpeMst} satisfy the modified restricted isometry property \eqref{SpeRIP}.
Again, without going into details,
we point out that the proof  is in the spirit of \cite{FouLai}:
for a fixed $\vZ \in \bR^{n \times n}$, 
establish a concentration inequality for $\left\|  (\va_i^\top \vZ \va_i)_{i=1}^m \right\|_1$ around its expectation $\sslash \vZ \sslash$,
 prove that this slanted norm is equivalent to the Frobenius norm, 
 and conclude with a covering argument.

\paragraph{Convergence of the idealized iterative hard thresholding.}
We now establish that the naive (and impractical) iterative hard thresholding algorithm \eqref{SimpleIHT} allows for stable and robust recovery of jointly low-rank and bisparse matrices under the standard restricted isometry property.
The precise statement appears after the important observation below.

\blem
\label{LemSimpleIHT}
Suppose that $\bm{\cA}: \bR^{n \times n} \to \bR^m$ satisfies the restricted isometry property \eqref{GenRIP} on $\Sigma_{(2s)}^{[2r]}$ with constant $\delta \in (0,1)$.
Then, for all $\vZ,\vZ' \in \Sigma_{(s)}^{[r]}$, one has
\be
\left| \langle \vZ, ( \bm{\cA}^* \bm{\cA} - \vI )(\vZ') \rangle \right| \le \delta  \|\vZ\|_F \|\vZ'\|_F.
\ee
\elem

\bpf
Assuming without loss of generality that $\|\vZ\|_F = \|\vZ'\|_F = 1$, we use in particular the parallelogram identity to write
\begin{align}
\left| \langle \vZ, ( \bm{\cA}^* \bm{\cA} - \vI)(\vZ') \rangle \right|
& = \left| \langle \bm{\cA}(\vZ), \bm{\cA}(\vZ') \rangle - \langle \vZ, \vZ' \rangle \right|\\
\nonumber
& = \left| \f{1}{4} \left( \|\bm{\cA}(\vZ + \vZ')\|_2^2 - \|\bm{\cA}(\vZ - \vZ')\|_2^2 \right) - \f{1}{4} \left( \|\vZ + \vZ'\|_F^2 - \| \vZ - \vZ'\|_F^2  \right) \right|\\
\nonumber
& \le \f{1}{4} \left| \|\bm{\cA}(\vZ + \vZ')\|_2^2 - \|\vZ + \vZ'\|_F^2  \right|
+  \f{1}{4} \left| \|\bm{\cA}(\vZ - \vZ')\|_2^2 - \|\vZ - \vZ'\|_F^2  \right|\\
\nonumber
& \le \f{1}{4} \delta  \|\vZ + \vZ'\|_F^2
+ \f{1}{4} \delta \|\vZ - \vZ'\|_F^2
= \f{1}{4} \delta  \left( 2 \|\vZ\|_F^2 + 2 \|\vZ'\|_F^2 \right) = \delta ,
\end{align}
which is the required result.
\epf

\bthm
\label{ThmSimpleIHT}
If the restricted isometry property \eqref{GenRIP} holds on $\Sigma_{(4s)}^{[4r]}$ with constant $\delta \in (0,1/2)$,
then any $\vX \in \Sigma_{(s)}^{[r]}$ is approximated from $\vy = \bm{\cA} \vX + \ve \in \bR^m$
as a cluster point $\vX_\infty$ of the sequence $(\vX_k)_{k\ge 0}$ defined by
\be
\vX_{k+1} = P_{(s)}^{[r]} \left( \vX_k + \bm{\cA}^*(\vy - \bm{\cA} \vX_k) \right)
\ee
with error
\be
\label{ObjThmSimpleIHT}
\|\vX - \vX_\infty\|_F \le C \|\ve\|_2 .
\ee
\ethm

\bpf
It is enough to prove that,
for all $k \ge 0$,
\be
\label{Obj1}
\|\vX - \vX_{k+1} \|_F \le \rho \|\vX - \vX_k \|_F + \tau \|\ve\|_2,
\qquad \mbox{with} \quad \rho:= 2 \delta< 1 \mbox{ and } \tau>0.
\ee
To start, notice that $\vX_{k+1}$ better approximates $\vX_k + \bm{\cA}^*(\vy - \bm{\cA} \vX_k) = \vX_k + \bm{\cA}^* \bm{\cA} (\vX - \vX_k)+ \bm{\cA}^*\ve$ 
as an element from $\Sigma_{(s)}^{[r]}$
than $\vX$ does,
so that
\be
\| \vX_k + \bm{\cA}^* \bm{\cA} (\vX - \vX_k) +\bm{\cA}^*\ve - \vX_{k+1} \|_F^2
\le 
\| \vX_k + \bm{\cA}^* \bm{\cA} (\vX - \vX_k) + \bm{\cA}^*\ve - \vX \|_F^2.
\ee
Introducing $\vX$ in the left-hand side,
expanding the squares, and simplifying leads to
\be
\label{Step1}
\|\vX - \vX_{k+1} \|_F^2 \le -2 \langle \vX - \vX_{k+1}, (\bm{\cA}^* \bm{\cA} - \vI )(\vX - \vX_k) + \bm{\cA}^*\ve\rangle.
\ee
Thanks to Lemma \ref{LemSimpleIHT}, we have
\be
\label{Step2}
| \langle \vX - \vX_{k+1}, (\bm{\cA}^* \bm{\cA} - \vI )(\vX - \vX_k) \rangle |
\le  2 \delta \| \vX - \vX_{k+1}\|_F  \|\vX - \vX_k \|_F,
\ee
while the restricted isometry property \eqref{GenRIP} also guarantees that 
\be
\label{Step3}
| \langle \vX - \vX_{k+1}, \bm{\cA}^*\ve\rangle |
= | \langle \bm{\cA} (\vX - \vX_{k+1}), \ve\rangle |
\le \|\bm{\cA} (\vX - \vX_{k+1})\|_2 \|\ve\|_2 
\le \sqrt{1+\delta} \| \vX - \vX_{k+1}\|_F \|\ve\|_2. 
\ee
Therefore, using \eqref{Step2} and \eqref{Step3} in \eqref{Step1}, we obtain
\be
\|\vX - \vX_{k+1} \|_F^2 \le 2 \delta \| \vX - \vX_{k+1}\|_F  \|\vX - \vX_k \|_F
+ \sqrt{1+\delta} \| \vX - \vX_{k+1}\|_F \|\ve\|_2,
\ee
which clearly implies the required estimates \eqref{Obj1} with $\tau = \sqrt{1+\delta}$
and \eqref{ObjThmSimpleIHT} with $C=\tau/(1-\rho)$.
\epf

\paragraph{The exact projection for $\Sigma_{(s)}^{[r]}$.}
Here, we prove the statement \eqref{FormProj} about the form of~$P_{(s)}^{[r]}$ 
before justifying that property \eqref{ProsT2}
holds for $T=P_{(s)}^{[r]}$.

\bprop
\label{LemFormProj}
For $\vM \in \bR^{n \times n}$,
the projection $P_{(s)}^{[r]}(\vM)$ of $\vM$ onto $\Sigma_{(s)}^{[r]}$
has the form $P^{[r]}(\vM_{S_\star \times S_\star})$,
where $S_\star$ maximizes $\|P^{[r]}(\vM_{S \times S})\|_F$
over all index sets $S$ of size $s$.
\eprop

\bpf
Let us remark that,
for any index set $T$,
\begin{align}
\label{Calc}
\|\vM - P^{[r]}(\vM_{T \times T}) \|_F^2
& = \|\vM_{\ol{T \times T}} + \vM_{T \times T} - P^{[r]}(\vM_{T \times T}) \|_F^2\\
\nonumber
& = \|\vM_{\ol{T \times T}}\|_F^2 + \|\vM_{T \times T} - P^{[r]}(\vM_{T \times T}) \|_F^2\\
\nonumber
& = \|\vM_{\ol{T \times T}}\|_F^2 + \|\vM_{T \times T}\|_F^2 - \| P^{[r]}(\vM_{T \times T}) \|_F^2\\
\nonumber
& = \|\vM \|_F^2  - \| P^{[r]}(\vM_{T \times T}) \|_F^2.
\end{align}
Now let $\vZ \in \Sigma_{(s)}^{[r]}$ and consider an index set $S$ of size $s$ such that $\vZ$ is supported on $S \times S$.
The defining property of $S_\star$, together with \eqref{Calc},
implies that
\begin{align}
\|\vM - P^{[r]}(\vM_{S_\star \times S_\star}) \|_F^2
& \le \|\vM \|_F^2  - \| P^{[r]}(\vM_{S \times S}) \|_F^2
 = \|\vM_{\ol{S \times S}}\|_F^2 + \|\vM_{S \times S} - P^{[r]}(\vM_{S \times S}) \|_F^2\\
 \nonumber
 & \le \|\vM_{\ol{S \times S}}\|_F^2 + \|\vM_{S \times S} - \vZ \|_F^2
 = \|\vM - \vZ \|_F^2,
\end{align}
where we have taken into account the facts that $P^{[r]}(\vM_{S \times S})$ is the best $r$-rank approximation to $\vM_{S \times S}$ and that $\vM_{\ol{S \times S}}$ and $\vM_{S \times S} - \vZ$ are disjointly supported.
Thus, we have proved that
$\|\vM - P^{[r]}(\vM_{S_\star \times S_\star})\|_F\le \|\vM - \vZ\|_F $
for all $\vZ \in \Sigma_{(s)}^{[r]}$, which is the desired result.
\epf

\bprop
\label{LemProjBiggerS}
For $\vM \in \bR^{n \times n}$,
considering an  index set $S_\star$ of size~$s$ with $P_{(s)}^{[r]}(\vM) = P^{[r]}(\vM_{S_\star \times S_\star})$, one has
\be
P_{(s)}^{[r]}(\vM) = P_{(s)}^{[r]}(\vM_{S' \times S'})
\qquad \mbox{whenever } S' \supseteq S_\star.
\ee
\eprop

\bpf
According to Proposition \ref{LemFormProj}, it is enough to verify that,
for any index set $S$ of size $s$,
\be
\left\| P^{[r]}( (\vM_{S' \times S'})_{S_\star \times S_\star} ) \right\|_F
\ge
\left\| P^{[r]}( (\vM_{S' \times S'})_{S \times S} ) \right\|_F.
\ee
But this is true because $(\vM_{S' \times S'})_{S_\star \times S_\star} = \vM_{S_\star \times S_\star}$
and 
$(\vM_{S' \times S'})_{S \times S} = (\vM_{S \times S})_{S' \times S'} $,
so that
\be
\left\| P^{[r]}( (\vM_{S' \times S'})_{S \times S} ) \right\|_F
\le \|P^{[r]}(\vM_{S \times S}) \|_F 
\le \|P^{[r]}(\vM_{S_\star \times S_\star})\|_F,
\ee
where the last inequality follows from the defining property of $S_\star$.
\epf


\vspace{-7mm}
\begin{thebibliography}{99}\vspace{-5mm}

\bibitem{Soh16}
S. Bahmani, and J. Romberg.
{\em Near-optimal estimation of simultaneously sparse and low-rank
  matrices from nested linear measurements.}
Information and Inference: A Journal of the IMA 5.3 (2016):331--351.

\bibitem{Blu}
T. Blumensath.
{\em Sampling and reconstructing signals from a union of linear subspaces.}
IEEE Transactions on Information Theory 57.7 (2011): 4660--4671.

\bibitem{CaiZha}
T. Cai and A. Zhang. 
{\em ROP: Matrix recovery via rank-one projections.} 
The Annals of Statistics 43.1 (2015): 102--138.

\bibitem{CanPla}
E. J. Cand\`es and Y. Plan. 
{\em Tight oracle inequalities for low-rank matrix recovery from a minimal number of noisy random measurements.}
IEEE Transactions on Information Theory 57.4 (2011): 2342--2359.

\bibitem{Can13}
E. J. Cand\`es, T. Strohmer, and V. Voroninski.
{\em Phaselift: exact and stable signal recovery from magnitude measurements via convex programming.} Communications on Pure and Applied Mathematics 66.8 (2013): 1241--1274.

\bibitem{DPW11}
R. DeVore, P. Guergana, and P. Wojtaszczyk, {\em Approximation of functions of few variables in high dimensions.} Constructive Approximation 33.1 (2011): 125-143.

\bibitem{For18}
  M. Fornasier, J. Maly, and V. Naumova.
  {\em Sparse PCA from inaccurate and incomplete measurements.}
  arXiv preprint, arXiv:1801.06240 (2018).

\bibitem{Fou19}
S. Foucart. {\em Sampling Schemes and Recovery Algorithms for Functions of Few Coordinate Variables.} preprint (2019).

\bibitem{FouLai}
S. Foucart and M.-J. Lai. 
{\em Sparse recovery with pre-Gaussian random matrices.} 
Studia Math 200.1 (2010): 91--102.

\bibitem{FouLec}
S. Foucart and G. Lecu\'e. 
{\em An IHT algorithm for sparse recovery from sub-exponential measurements.} 
IEEE Signal Processing Letters 24.9 (2017): 1280--1283.

\bibitem{fora13} 
S.~{F}oucart and H.~{R}auhut. 
{\em {A} {M}athematical {I}ntroduction to {C}ompressive {S}ensing}. 
{A}pplied and {N}umerical {H}armonic {A}nalysis. {B}irkh{\"a}user, 2013. 

\bibitem{FouSub}
S. Foucart and  S. Subramanian.
{\em Iterative hard thresholding for low-rank recovery from rank-one projections.}
Linear Algebra and its Applications 572 (2019): 117--134. 


\bibitem{Gep19}
J. Geppert, F. Krahmer, and D. St\"oger.
{\em Sparse power factorization: balancing peakiness and sample complexity.}
Advances in Computational Mathematics 45.3 (2019): 1711--1728.

\bibitem{GolDav}
M. Golbabaee and M. E. Davies. 
{\em Inexact gradient projection and fast data driven compressed sensing.} 
IEEE Transactions on Information Theory 64.10 (2018): 6707--6721.

\bibitem{HegIndSch}
C. Hegde, P. Indyk, and L. Schmidt. 
{\em Approximation algorithms for model-based compressive sensing.} 
IEEE Transactions on Information Theory 61.9 (2015): 5129--5147.

\bibitem{IwenVW}
M. Iwen, A. Viswanathan, and Y. Wang. 
{\em Robust sparse phase retrieval made easy.} 
Applied and Computational Harmonic Analysis 42.1 (2017): 135--142.

\bibitem{kakurate16} 
M.~{K}abanava, R.~{K}ueng, H.~{R}auhut, and U.~{T}erstiege. 
{\em {S}table low-rank matrix recovery via null space properties.} 
{I}nformation and {I}nference 5.4 (2016): 405--441. 

\bibitem{kurate17} 
R.~{K}ueng, H.~{R}auhut, and U.~{T}erstiege. 
{\em {L}ow rank matrix recovery from rank one measurements.} 
Applied and Computational Harmonic Analysis 42.1 (2017): 88--116. 

\bibitem{Lee17}
  K. Lee, Y. Wu, and Y. Bresler.
  {\em Near-optimal compressed sensing of a class of sparse low-rank matrices via sparse power factorization.}
  IEEE Transactions on Information Theory 64.3 (2017): 1666--1698.
  
\bibitem{MagdonIsmail} 
M. Magdon-Ismail,
{\em NP-hardness and inapproximability of sparse PCA.}
Information Processing Letters 126 (2017): 35--38. 

\bibitem{Manurangsi}
P. Manurangsi. 
{\em Almost-polynomial ratio ETH-hardness of approximating densest k-subgraph.}
STOC'17--Proceedings of the 49th Annual ACM SIGACT Symposium on Theory of Computing, ACM, pp. 954--961, 2017.

\bibitem{OJFEH}
S. Oymak, A. Jalali, M. Fazel, Y. C. Eldar, and B. Hassibi. 
{\em Simultaneously structured models with application to sparse and low-rank matrices.} 
IEEE Transactions on Information Theory 61.5 (2015): 2886--2908.

\bibitem{rascst17}
H. Rauhut, R. Schneider, and Z. Stojanac.
{\em Low rank tensor recovery via iterative hard thresholding.}
Linear Algebra and Its Applications 523 (2018): 220--262.

\bibitem{roklwuei16}
I. Roth, M. Kliesch, G. Wunder, and J. Eisert.
{\em Reliable recovery of hierarchically sparse signals and application in machine-type communications.}
arXiv preprint arXiv:1612.07806 (2016).

\bibitem{roflkueiwu18}
 I. Roth, A. Flinth, R. Kueng, J. Eisert, and G. Wunder.
{\em Hierarchical restricted isometry property for Kronecker product measurements.}
56th Annual Allerton Conference on Communication, Control, and Computing,
IEEE, pp. 632--638, 2018.
  

\bibitem{ShenLi}
J. Shen and P. Li.
{\em A tight bound of hard thresholding.} 
arXiv preprint arXiv:1605.01656 (2016).

\end{thebibliography}
\end{document}